\documentclass[reqno,11pt]{amsart}
\usepackage{amsmath}
\usepackage{amsfonts,amssymb}

\usepackage{latexsym}
\usepackage{amsfonts}
\usepackage{amsmath}
\usepackage{amssymb,mathrsfs}
\usepackage{bbm}
\usepackage[all]{xy}
\font\bbfont=msbm10 at 10pt
\font\sbbfont=msbm10 at 7pt
\def\CC{\mathscr C}
\def\C{\mbox{\bbfont C}}
\def\R{\mbox{\bbfont R}}
\def\1{\mathbbm 1}
\def\P{\mathbbm P}
\def\E{\mathbbm E}
\def\H{\mathcal H}
\def\HH{\mathscr H}
\def\Hc{{\mathscr H}_{\mbox{\sbbfont C}}}

\def\B{\mathcal B}
\def\F{\mathcal F}
\def\N{\mbox{\bbfont N}}
\def\eps{\sigma}
\def\e{\epsilon}
\def\tensor{\otimes}
\def\Gamm{\mathrm \Gamma}
\def\vac{\mathrm \Omega}
\def\te{\tau_\eps}
\def\z{\hat z}
\def\Cs{\CC(I,\eps)}
\def\Cd{\CC(I\times\{0,1\},\eps)}
\newtheorem{theorem}{Theorem}[section]
\newtheorem{proposition}[theorem]{Proposition}

\newtheorem{corollary}[theorem]{Corollary}

\newtheorem{lemma}[theorem]{Lemma}
\numberwithin{equation}{section}

\begin{document}

\title{Hypercontractivity in non-commutative holomorphic spaces}

\author[T. Kemp]{Todd Kemp}
\address{Cornell University \\ Malott Hall, Ithaca, NY \\ 14583-4201}
\email{tkemp@math.cornell.edu}
\date{\today}

\begin{abstract}
We prove an analog of Janson's strong hypercontractivity inequality in a class of non-commutative ``holomorphic" algebras. Our setting is the $q$-Gaussian algebras $\Gamm_q$ associated to the $q$-Fock spaces of Bozejko, K\"ummerer and Speicher, for $q\in[-1,1]$.  We construct subalgebras $\H_q\subset\Gamm_q$, a $q$-Segal-Bargmann transform, and prove Janson's strong hypercontractivity $L^2(\H_q)\to L^r(\H_q)$ for $r$ an even integer.
\end{abstract}

\maketitle

\section{Introduction}

As part of the work in the 1960s and 1970s to construct a mathematically consistent theory of interacting quantum fields,
Nelson proved his famous hypercontractivity inequality in its initial form \cite{Nelson 1}; by 1973 it evolved into the
following statement, which may be found in \cite{Nelson 2}.

\begin{theorem}[Nelson, 1973]\label{Nelson hyp thm} Let $A_\gamma$ be the Dirichlet form operator for Gauss measure
$d\gamma(x) = (2\pi)^{-n/2}e^{-|x|^2/2}dx$ on $\R^n$.  For $1< p\le r < \infty$ and $f\in L^p(\R^n,\gamma)$,
\begin{equation}\label{Nelson hyp}
\| e^{-tA_\gamma} f\|_r \le \|f\|_p,\quad\text{for}\quad t\ge t_N(p,r) = \frac{1}{2}\log\frac{r-1}{p-1}.
\end{equation}
For $t < t_N(p,r)$, $e^{-tA_\gamma}$ is not bounded from $L^p$ to $L^r$.
\end{theorem}
(If $p<2$, one must first extend $e^{-tA_\gamma}$ to $L^p$; this can be done, uniquely, and Theorem \ref{Nelson hyp thm} should
be interpreted as such in this case.  The same comment applies to all of the following.) It is worth noting that $t_N$,
the least time to contraction, does not depend on the dimension $n$ of the underlying space $\R^n$.

\medskip

The initial purpose of such hypercontractive inequalities was to prove the semiboundedness of Hamiltonians in the theory of Boson quantum fields.  (See, for example, \cite{Glimm}, \cite{Nelson 1}, and \cite{Segal 2}.)  In \cite{Gross 1}, Gross used this inequality (through an appropriate cut-off approximation) to show that the Boson energy operator in a model of 2-dimensional Euclidean quantum field theory has a unique ground state.  In that paper he also showed that if one represents the Fock space for Fermions as the $L^2$-space of a Clifford algebra (as in \cite{Segal 1}), then inequalities similar to \ref{Nelson hyp} also hold.  He developed this further in \cite{Gross 3}.

\medskip

Over the subsequent three decades, Nelson's hypercontractivity inequality (and its equivalent form, the logarithmic Sobolev inequality, invented by Gross in \cite{Gross 2}) found myriad applications in analysis, probability theory, differential geometry, statistical mechanics, and other areas of mathematics
and physics.  See, for example, the recent survey \cite{Gross 5}.

\medskip

The Fermion hypercontractivity inequality in \cite{Gross 3} remained unproven in its sharp form until the early 1990s.  Lindsay \cite{Lindsay} and Meyer \cite{Lindsay Meyer} proved that it holds $L^2\to L^r$ for $r=2,4,6,\ldots$ (and in the dual cases $L^{r'}\to L^2$ as well).  Soon after, Carlen and Lieb \cite{Carlen Lieb} were able to complete Gross' original argument with some clever non-commutative integration inequalities, thus proving the full result.  (Precisely: they showed that the Clifford algebra analogs of the inequalities \ref{Nelson hyp} hold with exactly the same constants.)

\medskip

Then, in 1997, Biane \cite{Biane 1} extended Carlen and Lieb's work beyond the Fermionic (Clifford algebra) setting to the $q$-Gaussian von Neumann algebras $\Gamm_q$ of Bozejko, K\"ummerer, and Speicher \cite{BKS}.  His theorem may be stated as follows.
\begin{theorem}[Biane, 1997]\label{Free hyp thm} Let $-1<q<1$, let $N_q$ denote the number operator associated to $\Gamm_q$, and let
$\|\cdot\|_p$ be the non-commutative $L^p$-norm associated to the vacuum expectation state $\tau_q$ on $\Gamm_q$.  Then for $1<p\le r <\infty$,
\[ \| e^{-tN_q}f \|_r \le \|f\|_p\text{ for all }f\in L^p(\Gamm_q,\tau_q) \quad \text{iff} \quad t \ge t_N(p,r). \]
\end{theorem}
Of particular interest is the case $q=0$ which corresponds to free probabiliy.  Biane proved the full result (for $-1<q<1$)
by first extending Carlen and Lieb's work to the case of a system of mixed spins (in a von Neumann algebra generated
by elements which satisfy some commutation and some anti-commutation relations), and then applying a central limit
theorem due to Speicher \cite{Speicher}.  The case $q=-1$ is Carlen and Lieb's adaptation of Gross' work, while the $q=1$
case is Nelson's original hypercontractive estimate (Theorem \ref{Nelson hyp thm}).

\medskip

Concurrent to the work on non-commutative hypercontractivity, a different sort of extension of Nelson's theorem was being developed.  In 1983 Janson, \cite{Janson}, discovered that if one restricts the semigroup $e^{-tA_\gamma}$ in Theorem \ref{Nelson hyp thm} to holomorphic functions on $\R^{2n} \cong \C^n$ then the contractivity of Equation \ref{Nelson hyp} is attained in a shorter time than $t_N$.  Writing $\H L^p = L^p(\R^{2n},\gamma)\cap\text{Hol}(\C^n)$, Janson's  {\it strong hypercontractivity} may be stated thus.
\begin{theorem}[Janson, 1983]\label{Janson hyp thm} Let $0< p \le r < \infty$, and let $f \in \H L^p$.  Then
\begin{equation}\label{Janson hyp}
\| e^{-tA_\gamma} f \|_r \le \|f\|_p,\quad\text{for}\quad t\ge t_J(p,r)=\frac{1}{2}\log\frac{r}{p}.
\end{equation}
For $t < t_J(p,r)$, $e^{-tA_\gamma}$ is not bounded from $\H L^p$ to $\H L^r$.
\end{theorem}
Note that the least time $t_J$ to contraction is shorter than the time $t_N$ (if $1 < p < r < \infty$).  Moreover, Janson's result
holds as $p\to 0$, in a regime where the semigroup $e^{-tA_\gamma}$ is not even well-defined in the full $L^p$-space.
These results have been further generalized by Gross in \cite{Gross 4} to the case of complex manifolds.

\medskip

In this paper, non-commutative algebras $\H_q$ will be introduced, which are $q$-deformations of the algebra of holomorphic functions.  The special cases $q=\pm 1$ and $q=0$ are already known; $\H_{-1}$ is defined in \cite{BSZ}, while $\H_0$ is isomorphic to the free Segal-Bargmann space of \cite{Biane 2}.  We will construct a unitary isomorphism ${\mathscr S}_q$ from $L^2(\Gamm_q)$ to $L^2(\H_q)$, which is a $q$-analog of the Segal-Bargmann transform.  $\H_q$ itself will be constructed as a subalgebra of $\Gamm_q$, and so inherits its $p$-norms as well as its number operator $N_q$.  In the context of these $q$-deformed Segal Bargmann spaces, the following theorem is our main result.
\begin{theorem}\label{Main Theorem} For $-1 \le q < 1$ and $r$ an even integer,
\[ \| e^{-tN_q}f \|_r \le \|f\|_2\text{ for all }f\in L^2(\H_q,\tau_q) \quad \text{iff} \quad t \ge t_J(2,r). \]
\end{theorem}
It is interesting that the least time to contraction, $t_J$, is independent of both the dimension of the underlying space and the parameter $q$. We fully expect the same results to hold $L^p(\H_q)\to L^r(\H_q)$ for $2 \le p \le r < \infty$, but standard interpolation techniques fail to work in the holomorphic algebras we consider.  (In particular, the dual results that Lindsay and Meyer achieved in the full Clifford algebra do not follow in this holomorphic setting.) 
\medskip

This paper is organized as follows.  We begin with a summary of the $q$-Fock spaces $\F_q$ and the von Neumann algebras
$\Gamm_q$ associated to them.  We will also define the holomorphic subalgebras $\H_q$ and construct a $q$-Segal-Bargmann transform.  In the subsequent section, we prove the appropriate strong hypercontractivity estimates for algebras with arbitrary mixed spins (mixed commutation and anti-commuations relations), much in the spirit of Biane's approach \cite{Biane 1}.  We then proceed to review Speicher's central limit theorem, and apply it to prove Theorem \ref{Main Theorem}.

\section{The $q$-Fock space and associated algebras}\label{q-Algebras}

We begin by briefly reviewing the $q$-Fock spaces of Bozejko, K\"ummerer and Speicher, relevant aspects of the von Neumann algebras $\Gamm_q$ (which are related to the creation and annihilation operators on $\F_q$), and the number operators on them.  We then proceed to define the Banach algebra $\H_q$ which corresponds to the classical Segal-Bargmann space, and exhibit a $\ast$-isomorphism between $\H_0$ and the free Segal-Bargmann space ${\mathscr C}_{hol}$ defined in \cite{Biane 2}.  We finally construct a generalized $q$-Segal-Bargmann transform, which is a unitary isomorphism $L^2(\Gamm_q)\to L^2(\H_q)$ that respects the action of the number operator.

\subsection{The $q$-Fock space $\F_q$ and the algebra $\Gamm_q$}\label{q-Fock}
Our development closely follows that found in \cite{Biane 1}; the details may be found in \cite{BKS}.  Let $\HH$ be a {\it real} Hilbert space with complexification $\Hc$. Let $\vac$ be a unit vector in a $1$-dimensional complex Hilbert space (disjoint from $\Hc$).  We refer to  $\vac$ as \emph{the vacuum}, and by convention define $\Hc^{\tensor 0} \equiv \C\vac$.  The \emph{algebraic Fock space} $\F(\HH)$ is defined as
\[ \F(\HH) \equiv \bigoplus_{n=0}^\infty \Hc^{\tensor n}, \]
where the direct sum and tensor product are algebraic.  For any $q\in[-1,1]$, we then define a Hermitian form $(\cdot,\cdot)_q$ to be the conjugate-linear extension of
\begin{eqnarray*}
                                                                             (\vac,\vac)_q & = & 1 \\
(f_1\tensor\cdots\tensor f_j,g_1\tensor\cdots\tensor g_k)_q & = & \delta_{jk}\sum_{\pi\in\mathcal{S}_k}q^{\iota(\pi)}(f_1,g_{\pi 1})\cdots(f_k,g_{\pi k}),
\end{eqnarray*}
for $f_i,g_i\in\Hc$, where $\mathcal{S}_k$ is the symmetric group on $k$ symbols, and $\iota(\pi)$ counts the number of inversions in $\pi$; that is
\[ \iota(\pi) = \#\{(i,j)\,;\,1\le i < j \le k, \pi i > \pi j\}. \]
The reader may readily verify that $(-1)^{\iota(\pi)} = \text{parity}(\pi)$ for any permutation $\pi$.  Hence, the form $(\cdot,\cdot)_{-1}$ reduces
to the standard Hermitian form associated to the Fermion Fock space.  Similarly, the form $(\cdot,\cdot)_{1}$ yields the standard Hermitian form on the Boson Fock space.  In each of these cases the form is degenerate, thus requiring that we take a quotient of $\F(\HH)$ before completing to form the Fermion or Boson Fock space.  It is somewhat remarkable that, for $-1<q<1$, the form $(\cdot,\cdot)_q$ is already non-degenerate on $\F(\HH)$.
\begin{proposition}[\cite{BKS}]\label{q-form is nondegenerate}
The Hermitian form $(\cdot,\cdot)_q$ is positive semi-definite on $\F(\HH)$.  Moreover, it is an inner product on $\F(\HH)$ for $-1<q<1$.
\end{proposition}
For $-1<q<1$, the \emph{$q$-Fock space} $\F_q(\HH)$ is defined as the completion of $\F(\HH)$ with respect to the inner-product $(\cdot,\cdot)_q$. (It should be noted that, in the case $q=0$, the definition of the form $(\cdot,\cdot)_0$ requires the convention that $0^0 = 1$.  It follows that $\F_0(\HH)$ is just $\bigoplus_{n=0}^\infty\Hc^{\tensor n}$ with the Hilbert space tensor product and direct sum.)  These spaces interpolate between the classical Boson and Fermion Fock spaces $\F_{\pm 1}(\HH)$, which are constructed by first taking the quotient of $\F(\HH)$ by the kernel of $(\cdot,\cdot)_{\pm 1}$ and then completing.

\medskip

As in the classical theory, the spaces $\F_q$ come equipped with creation and annihilation operators. For any vector $f\in\HH\subset\Hc$, define the \emph{creation operator} $c_q(f)$ on $\F_q(\HH)$ to extend
\begin{eqnarray*}
                                 c_q(f)\vac & = & f \\
c_q(f)f_1\tensor\cdots\tensor f_k & = & f\tensor f_1\tensor\cdots\tensor f_k.
\end{eqnarray*}
The \emph{annihilation operator} $c_q^\ast(f)$ is its adjoint, which the reader may compute satisfies
\begin{eqnarray*}
                                c_q^\ast(f)\vac & = & 0 \\
c_q^\ast(f)f_1\tensor\cdots\tensor f_k & = & \sum_{j=1}^k q^{j-1}(f_j,f)f_1\tensor\cdots\tensor f_{j-1}\tensor f_{j+1}\tensor\cdots\tensor f_k.
\end{eqnarray*}
These are similar to the definitions of the creation and annihilation operators in the Fermion and Boson cases, where appropriate (anti)symmetrization must also be applied.  One notable difference is that, in the Boson ($q=1$) case, the operators are unbounded. For $q<1$, the creation and annihilation operators are always bounded, and hence we may discuss the von Neumann algebra they generate without difficulties.

\medskip

The operators $c_q,c_q^\ast$ satisfy the \emph{$q$-commutation relations}, which interpolate between the canonical commutation relations (CCR)
and canonical anticommutation relations (CAR) usually associated to the Boson and Fermion Fock spaces.  Over the $q$-Fock space, we have
\begin{equation}\label{q-commutation relations}
c_q^\ast(g)c_q(f) - qc_q(f)c_q^\ast(g) = (f,g)\text{id}_{\F_q(\HH)}\;\;\text{for}\;\; f,g\in\HH.
\end{equation}
It is worth pausing at this point to note one significant difference between the $q=\pm 1$ cases and the $-1<q<1$ cases.  For both Bosons and
Fermions, the operators $c,c^\ast$ also satisfy additional (anti)commutation relations.  In the Boson case, for example, $c(f)$ and $c(g)$ commute for any choices of $f$ and $g$.  It is a fact, however, that if $q\ne\pm 1$ there are {\it no relations} between $c_q(f)$ and $c_q(g)$ if $(f,g)=0$.

\medskip

It is a well-known theorem that the creation and annihilation operators in the Boson and Fermion cases are irreducible; that is, they have no non-trivial invariant subspaces. That theorem is also true for the operators $c_q$ for $-1<q<1$, although a published proof does not seem to exist.  We prove it here for completeness.  The $q=0$ case will be used in Proposition \ref{* isom} below.
\begin{theorem}\label{irreducible}
For $-1<q<1$, the von Neumann algebra generated by $\{c_q(h)\,;\,h\in\HH\}$ is $\B(\F_q(\HH))$.
\end{theorem}
\begin{proof} Denote by ${\mathscr W}_q$ the von Neumann algebra generated by the $c_q$'s. We consider the $q=0$ case first.  Let $\{e_1,e_2,\ldots\}$ be an
orthonormal basis for $\HH$, and consider the operator
\[P = \sum_{j=1}^\infty c(e_j)\,c^\ast(e_j) \]
(where $c(h) = c_0(h)$), which is in ${\mathscr W}_0$ since ${\mathscr W}_0$ is weakly closed. It is easy to calculate that $P(e_{i_1}\tensor\cdots\tensor e_{i_n}) = e_{i_1}\tensor\cdots\tensor e_{i_n}$, while $P\vac = 0$.  Thus, $P$ is the projection onto the orthogonal complement of the vacuum.  So ${\mathscr W}_0\ni 1-P = P_\vac$, the projection onto the vacuum.  Therefore ${\mathscr W}_0$ contains the operator
\[ c(e_{i_1})\cdots c(e_{i_n})P_{\vac}c^\ast(e_{j_1})\cdots c^\ast(e_{j_m}),  \]
which is the rank-$1$ operator with image spanned by $e_{i_1}\tensor\cdots\tensor e_{i_n}$ and kernel orthogonal to $e_{j_1}\tensor\cdots\tensor e_{j_m}$.
It follows that ${\mathscr W}_0$ contains all finite rank operators, and hence is
the full algebra $\B(\F_0(\HH))$.

\medskip

For $q\ne 0$, it is proved in \cite{DN} that there is a unitary map $U_q\colon\F_0\to\F_q$, which preserves the vacuum and satisfies
\[ U_q {\mathscr C}_0 U_q^\ast \subseteq {\mathscr C}_q, \]
where ${\mathscr C}_q$ is the $C^\ast$-algebra generated by $\{c_q(h)\,;\,h\in\HH\}$.  As ${\mathscr W}_q$ is the weak closure of ${\mathscr C}_q$, it follows easily that $\B(\F_q(\HH))=U_q\B(\F_0(\HH))U_q^\ast = U_q{\mathscr W}_0 U_q^\ast \subseteq {\mathscr W}_q$ as well, and this completes the proof.  \end{proof}

\medskip

For $q<1$ and for each $f\in\HH$, define the self-adjoint operator $X_q(f)$ on $\F_q(\HH)$ by $X_q(f) = c_q(f) + c_q^\ast(f)$.  These operators are in ${\mathscr W}_q = \B(\F_q(\HH))$, but they do not generate it.  The von Neumann algebra they do generate is defined to be $\Gamm_q(\HH)$, the \emph{$q$-Gaussian algebra} over $\HH$.  (In the $q=1$ case, $\Gamm_1(\HH)$ is the von Neumann algebra generated by the operators $\varphi(X(f))$ for $\varphi\in L^\infty(\R)$.)  The notation $\Gamm_q$ is chosen to be consistent with the {\it second quantization functor} from constructive quantum field theory (see \cite{BSZ}), which assigns to each real Hilbert space $\HH$ a von Neumann
algebra $\Gamm(\HH)$ and to each contraction $T\colon\HH\to\mathscr K$ a unital positivity-preserving map $\Gamm(T)\colon\Gamm(\HH)\to\Gamm(\mathscr K)$.  Indeed, $\Gamm_q$ can be construed as such a functor as well. 

\medskip

The isomorphism classes of the von Neumann algebras $\Gamm_q(\HH)$ for $q\notin\{\pm 1, 0\}$ are not yet understood.  (For some partial results, however, see \cite{Ricard} and \cite{Sniady}.)  The $\pm 1$ cases have been understood since antiquity: $\Gamm_1(\HH) = L^\infty(M,\gamma)$ for a certain measure space $M$ with a Gaussian measure $\gamma$ , while $\Gamm_{-1}(\HH)$ is a Clifford algebra modeled on $\HH$.  These facts rely upon the additional commutation relations between $c(f)$ and $c(g)$ that hold in those cases. (Indeed, in the Boson case $X(f)$ and $X(g)$ commute, resulting in a commutative von Neumann algebra $\Gamm(\HH)$. It is primarily for this reason that it is customary to begin with a real Hilbert space and complexify -- if $c(f)$ were defined for all $f\in\Hc$, then $c(f)$ and $c(g)$ would no longer commute even in the Boson case.  While there are no commutation relations between $c_q(f)$ and $c_q(g)$, it is still advantageous for us to have the real subspace $\HH\subset\Hc$ in order to define the holomorphic subalgebra in section \ref{q-holomorphic}.)  $\Gamm_0(\HH)$ was shown (in \cite{Voiculescu}) to be isomorphic to the group von Neumann algebra of a free group with countably many generators.
\medskip

One known fact about the algebras $\Gamm_q(\HH)$ for $-1<q<1$ is that they are all type $II_1$ factors.  This is a consequence (in the $\dim\HH = \infty$ case) of the following theorem, which was proved in \cite{BSp}. \begin{proposition}[Bozejko, Speicher]\label{trace prop} Let $-1<q<1$. The vacuum expectation state $\tau_q(A) = (A\vac,\vac)_q$ on $\B(\F_q(\HH))$ restricts to a faithful, normal, finite trace on $\Gamm_q(\HH)$.
\end{proposition}
The reader may wish to verify that $\tau_q(c_q^\ast c_q) = 1$, while $\tau_q(c_qc_q^\ast)=0$; hence, $\tau_q$ is certainly not a trace on all of $\B(\F_q(\HH))$.
\medskip

The algebra $\Gamm_q$ can actually be included as a dense subspace of $\F_q$.  The map $A\mapsto A\vac$ is one-to-one from $\Gamm_q$ into
$\F_q$.  The precise action of this map will be important to us, and so it bears mentioning.  The \emph{$q$-Hermite polynomials} $H^q_n$ are one-variable real polynomials defined so that $H^q_0(x) = 1$, $H^q_1(x) = x$, and satisfying the following recurrence relation:
\begin{equation}\label{q-Hermite}
xH^q_n(x) = H^q_{n+1}(x) +\frac{q^n-1}{q-1}H^q_{n-1}(x),
\end{equation}
where $(q^n-1)/(q-1)$ is to be interpreted as $n$ when $q=1$.  In this case, the generated polynomials $H^1_n$ are precisely the Hermite
polynomials that play an important role in the Boson theory.  When $q=0$, the polynomials $H^0_n$ are the Tchebyshev polynomials, and play
an analogous role in the theory of semi-circular systems (see \cite{Voiculescu}). We can express the action of the above map $A\mapsto A\vac$
succinctly in terms of the polynomials $H^q_n$.  The following proposition is proved in \cite{BKS}.
\begin{proposition}\label{Wick ordering 1}
The map $A\mapsto A\vac$ from $\Gamm_q$ to $\F_q$ is one-to-one, and extends to a unitary isomorphism $L^2(\Gamm_q,\tau_q)\to\F_q$.
If $\{e_j\}$ are orthonormal vectors in $\HH$ and $j_{\ell} \ne j_{\ell+1}$ for $1 \le \ell \le k-1$, then
\begin{equation}\label{q-Hermite mapping}
H^q_{n_1}(X_q(e_{j_1}))\cdots H^q_{n_k}(X_q(e_{j_k}))\vac = e_{j_1}^{\tensor n_1}\tensor\cdots\tensor e_{j_k}^{\tensor n_k}.
\end{equation}
\end{proposition}

\medskip

The algebraic Fock space $\F(\HH)$ carries a number operator $N$, whose action is given by
\begin{eqnarray*}
N\vac & = & 0 \\
N(f_1\tensor\cdots\tensor f_n) & = & nf_1\tensor\cdots\tensor f_n.
\end{eqnarray*}
This operator extends to a densely-defined, essentially self-adjoint operator $N_q$ on $\F_q(\HH)$.  The algebra $\Gamm_q$ then inherits
the action of $N_q$, via the map in Proposition \ref{Wick ordering 1}. The reader may readily check that if $\{e_j\}$ are orthonormal vectors in $\HH$ and $j_{\ell} \ne j_{\ell+1}$ for $1 \le \ell \le k-1$, then the element
\[ H^q_{n_1}(X_q(e_{j_1}))\cdots H^q_{n_k}(X_q(e_{j_k}))  \]
is an eigenvector of $N_q$ with eigenvalue $n_1+\cdots+n_k$.  In the case $\HH=\R^d$, this is a precise analogy to the action of the number operator for Bosons.  The algebra $\Gamm_1(\R^d)$ is isomorphic to $L^\infty(\R^d,\gamma)$, and the operators $X_1(e_j)$ (for the standard
basis vectors $e_j$) are multiplication by the coordinate functions $x_j$.  The Boson number operator $A_\gamma$ has $H^1_{n_1}(x_1)\cdots H^1_{n_k}(x_k)$ as an eigenvector, with eigenvalue $n_1+\cdots+n_k$.

\medskip 

The number operator generates a contraction semigroup $e^{-tN_q}$ on $L^2(\Gamm_q,\tau_q)$, which is known to restrict
for $p >2$, and extend for $1\le p <2$, to a contraction semigroup on $L^p(\Gamm_q,\tau_q)$.  Biane's hypercontractivity theorem, Theorem \ref{Free hyp thm}, is an extension of these results.

\subsection{the holomorphic algebra, and the $q$-Segal-Bargmann transform}\label{q-holomorphic}
Let $q<1$.  We wish to define a Banach algebra of ``holomorphic" elements in $\B(\F_q)$.  To that end, we follow a similar procedure to the formal construction of holomorphic polynomials.  We begin by doubling the number of variables, and so we consider the algebra $\Gamm_q(\HH\oplus\HH)$.  This algebra contains two independent copies of the variable $X(h)\in\Gamm_q(\HH)$: $X(h,0)$ and $X(0,h)$.  (Here, $(h,0)$ denotes a pair in $\HH\oplus\HH$, not the inner product of $h$ with $0$. Whenever this ambiguity in notation may be confusing, we will clarify by denoting the inner product as $(\cdot,\cdot)_{\mathscr K}$ for the appropriate Hilbert space $\mathscr K$.)  We then introduce a new variable $Z(h)$,
\begin{equation}\label{Z}
Z(h) = \frac{1}{\sqrt{2}}(X(h,0) + iX(0,h)).
\end{equation}
In the case $\HH=\R$ and $q=1$, this precisely corresponds to the holomorphic variable $z = (x+iy)/\sqrt{2}$, the normalization chosen
so that $z$ is a unit vector in $\H L^2(\gamma)$.  We define the {\bf $q$-holomorphic algebra} $\H_q(\Hc)$ as the Banach algebra generated by
$\{Z(h)\,;\,h\in\HH\}$.

\medskip

In \cite{Biane 2}, Biane introduced a Banach algebra ${\mathscr C}_{hol}$ in the $q=0$ case which is also an analog of the algebra of holomorphic
functions.  His algebra is not contained in $\Gamm_0(\HH\oplus\HH)$, so it is less natural to consider an action of $N_0$ on it.  We introduce it here (with slightly changed notation to avoid inconsistencies) to show that it is isomorphic to $\H_0$, and so the work presented here indeed generalizes Biane's results.  Consider the von Neumann algebra $\B(\F_0(\HH\oplus\HH))$; it contains all the operators $c_0(h,g)$ and their adjoints, for $h,g\in\HH$.  Define the
operator
\[ B(h) = c_0(h,0) + c_0^\ast(0,h). \]
Let ${\mathscr C}(\Hc)$ be the von Neumann algebra generated by $\{B(h)\,;\,h\in\HH\}$, and ${\mathscr C}_{hol}(\Hc)$ the Banach algebra so generated. The vacuum expectation state $\tau_0(A) = (A\vac,\vac)$ restricts to a faithful, normal, finite trace on ${\mathscr C}(\Hc)$, and the map $h\mapsto B(h)$ is a circular system with respect to $\tau_0$ (see \cite{Voiculescu}).

\medskip

Although Biane's algebra ${\mathscr C}_{hol}(\Hc)$ is not contained in $\Gamm_0(\HH\oplus\HH)$, it is in fact isomorphic to our algebra $\H_0(\Hc)$, in the following strong sense.
\begin{proposition}\label{* isom}
There is a $\ast$-automorphism of  $\B(\F_0(\HH\oplus\HH))$ which maps $\Gamm_0(\HH\oplus\HH)$ onto ${\mathscr C}(\Hc)$.  In particular, it maps $Z(h)$ to $B(h)$, and so sends $\H_0(\Hc)$ to ${\mathscr C}_{hol}(\Hc)$.
\end{proposition}

\begin{proof} By Theorem \ref{irreducible}, we can define an endomorphism of $\B(\F_0(\HH\oplus\HH))$ on the generators $c_0(h,g)$ by
\begin{eqnarray*}
\alpha(c_0(h,g)) & = & \frac{1}{\sqrt{2}}(c_0(h,h) + ic_0(-g,g)) \\
\alpha(c_0^\ast(h,g)) & = & \frac{1}{\sqrt{2}}(c_0^\ast(h,h) - i c_0^\ast(-g,g)).
\end{eqnarray*}
A straightforward computation verifies that the operators $\alpha(c_0(h,g))$ satisfy the $0$-commutation relations of Equation \ref{q-commutation relations}.
Hence, $\alpha$ extends to a $\ast$-homomorphism.  It can also easily be checked that $\alpha$ has an inverse of the form
\begin{eqnarray*}
\alpha^{-1}(c_0(h,g)) & = & \frac{1}{\sqrt{2}}(c_0(h+g,0)+ic_0(0,h-g)) \\
\alpha^{-1}(c_0^\ast(h,g)) & = & \frac{1}{\sqrt{2}}(c_0^\ast(h+g,0)-ic_0^\ast(0,h-g)),
\end{eqnarray*}
which also extends to a $\ast$-homomorphism.  Hence, $\alpha$ is a $\ast$-automorphism.  Finally, one can calculate that $\alpha(Z(h)) = B(h)$.
Whence, $\alpha$ maps $\H_0$ onto ${\mathscr C}_{hol}$, and so maps $W^\ast(\H_0(\Hc)) = \Gamm_0(\HH\oplus\HH)$ onto ${\mathscr C}(\Hc)$.
\end{proof}
It should be noted that Proposition \ref{* isom} generalizes automatically to $q\ne 0$; however, we are only concerned with the $q=0$ case for
${\mathscr C}_{hol}$.

\begin{corollary}\label{p-norms same}
The map $\alpha$ from $\Gamm_0(\HH\oplus\HH)$ to ${\mathscr C}(\Hc)$ extends to an isometric isomorphism $L^p(\Gamm_0,\tau_0)\to L^p({\mathscr C},\tau_0)$ for $0 < p \le \infty$.
\end{corollary}

\begin{proof} Since $\alpha$ is a $\ast$-automorphism of the full von Neumann algebra $\B(\F_0)$, by Wigner's theorem it is induced (through conjugation) by either a unitary or an anti-unitary on $\F_0$.  Suppose it is a unitary, $U$, so that $\alpha(A) = U^\ast A U$ for each $A\in\B(\F_0)$. Recall that $\tau_0$ is known to restrict to a tracial state on both $W^\ast(\H_0)$ and $\mathscr C$.  Hence, for $A\in\H_0$ and $p>0$,
\[ \|\alpha(A)\|_p^p = \tau_0(|\alpha(A)|^p) = \tau_0(\alpha(|A|^p)) = \tau_0(U^\ast|A|^pU) = \tau_0(|A|^p) = \|A\|_p^p. \]
It follows that $\alpha$ extends to an isometric isomorphism from $L^p(\H_0,\tau_0)$ onto $L^p({\mathscr C}_{hol},\tau_0)$.  The anti-unitary case is similar. \end{proof}
Hence, the algebraic map which sends $Z(h)$ to $B(h)$ preserves all $L^p$ topology (even for $p<1$), and so the analyses of the spaces $\H_0$ and ${\mathscr C}_{hol}$ are very much the same. 

\medskip

In the commutative context, one of the most powerful tools in this area is the \emph{Segal-Bargmann transform} $\mathscr S$, which is a unitary
isomorphism
\[  {\mathscr S}\colon L^2(\R^d,\gamma)\to\H L^2(\C^d,\gamma'). \]
Here, $\gamma'$ denotes the measure whose density with respect to Lebesgue measure is (a constant multiple of) the complex Gaussian $\exp(-|z|^2)$, rather than $\exp(-|x|^2/2)$ as in $\gamma$. The Hermite polynomials $H^1_{n_1}(x_1)\cdots H^1_{n_d}(x_d)$, appropriately normalized, form an orthonormal basis of $L^2(\R^d,\gamma)$, and the action of $\mathscr S$ on this basis is simple:
\begin{equation}\label{action of S}
{\mathscr S} \colon H^1_{n_1}(x_1)\cdots H^1_{n_d}(x_d) \mapsto z_1^{n_1}\cdots z_d^{n_d}.
\end{equation}
So $\mathscr S$ maps the Hermite polynomials to the holomorphic monomials.

\medskip

(A note on normalization.  Instead of changing the measure $\gamma\to\gamma'$, we could redefine
${\mathscr S}'\colon L^2(\R^d,\gamma)\to\H L^2(\C^d,\gamma)$ by setting ${\mathscr S}'f(z) = {\mathscr S}f(z/\sqrt{2})$.  This map is, of course,
a unitary isomorphism.  It is this point of view that we take while generalizing the Segal-Bargmann transform.  After all, we have already built the
factor $1/\sqrt{2}$ in to the variable $Z(h)$.)

\medskip

In \cite{Biane 2}, a free Segal-Bargmann transform is introduced, which is a unitary isomorphism $L^2(\Gamma_0,\tau_0)\to L^2({\mathscr C}_{hol},\tau_0)$.  We will modify this transform and extend it to all $q\in[-1,1]$, and further show that besides generalizing the classical
transform $\mathscr S$ it respects the action of the number operator.  First, we will need to understand the embedding of $\H_q(\Hc)$ in
$\F_q(\HH\oplus\HH)$ (it is injected via the map $A\mapsto A\vac$, which is one-to-one on all of $\Gamm_q(\HH\oplus\HH)$ by
Proposition \ref{Wick ordering 1}).  Consider the diagonal mapping $\delta\colon\Hc\to\Hc\oplus\Hc$ defined $\delta(h) = 2^{-1/2}(h,ih)$.  Since
$\delta$ is isometric, it extends to an isometric embedding $\delta_q\colon\F_q(\HH)\hookrightarrow\F_q(\HH\oplus\HH)$ (that is, $\delta_q(h_1\tensor h_2\tensor\cdots) = \delta(h_1)\tensor\delta(h_2)\tensor\cdots$).

\begin{proposition}\label{Wick ordering 2}
The map $A\mapsto A\vac$ injecting $\H_q(\Hc)\hookrightarrow\F_q(\HH\oplus\HH)$ extends to a unitary isomorphism
$L^2(\H_q(\Hc),\tau_q)\to\delta_q\F_q(\HH)$.  If $\{e_j\}$ are orthonormal vectors in $\HH$, then
\begin{equation}\label{monomial mapping}
Z_q(e_{j_1})^{n_1}\cdots Z_q(e_{j_k})^{n_k}\vac = \delta_q(e_{j_1}^{\tensor n_1}\tensor\cdots\tensor e_{j_k}^{\tensor n_k}).
\end{equation}
\end{proposition}

\begin{proof} Let $\phi = h_1\tensor\cdots\tensor h_n\in\F(\HH)$, and consider $Z_q(h)\delta(\phi)$.  We may compute
\begin{eqnarray*}
X_q(h,0)\delta(\phi) & = & (c_q(h,0)+c_q^\ast(h,0))\delta(\phi) \\
                            & = & (h,0)\tensor\delta(\phi) + \sum_{j=1}^n q^{j-1}\left(2^{-1/2}(h_j,ih_j),(h,0)\right)_{\HH\oplus\HH}\delta(\hat\phi_j) \\
                            & = & (h,0)\tensor\delta(\phi) + \frac{1}{\sqrt{2}}\sum_{j=1}^n q^{j-1}(h_j,h)_{\HH}\delta(\hat\phi_j), 
\end{eqnarray*}
where $\hat\phi_j = h_1\tensor\cdots\tensor h_{j-1}\tensor h_{j+1}\tensor\cdots\tensor h_n$.  A similar calculation shows that
\[
X_q(0,h)\delta(\phi) = (0,h)\tensor\delta(\phi) + \frac{1}{\sqrt{2}}\sum_{j=1}^n q^{j-1}(ih_j,h)_{\HH}\delta(\hat\phi_j),
\]
and so in the sum $X_q(h,0)+iX_q(0,h)$ the $c_q^\ast$ terms cancel.  (Note, we have assumed as is standard that the complexified inner product $(h,g)$
is linear in $h$ and conjugate-linear in $g$.)  Thus, we have $Z_q(h)\delta(\phi) = 2^{-1/2}(h,ih)\tensor\delta(\phi) = \delta(h\tensor\phi)$.
Equation \ref{monomial mapping} now follows by induction, and the theorem follows since such vectors are dense in $\delta_q\F_q(\HH)$. 
\end{proof}

We may now define the {\bf $q$-Segal-Bargmann transform} as follows.  Propositions \ref{Wick ordering 1} and \ref{Wick ordering 2} give (up to
the map $\delta_q$) unitary equivalences between the Fock space $\F_q(\HH)$ and both $L^2(\Gamm_q(\HH),\tau_q)$ and $L^2(\H_q(\Hc),\tau_q)$.
The $q$-Segal-Bargmann transform ${\mathscr S}_q$ is the composition of these unitary isomorphisms.  That is, ${\mathscr S}_q$ is the unitary map
which makes the following diagram commute.
\[
\xymatrix{
\F_q(\HH) \ar@{^{(}->}[r]^{\delta_q} & \F_q(\HH\oplus\HH) \\
L^2(\Gamm_q(\HH),\tau_q) \ar[u]^{A\mapsto A\vac} \ar[r]_{{\mathscr S}_q} & L^2(\H_q(\Hc),\tau_q) \ar[u]_{A\mapsto A\vac}
}
\]
By Equations \ref{q-Hermite mapping} and \ref{monomial mapping}, the action of ${\mathscr S}_q$ can be expressed in terms of the $q$-Hermite polynomials.  If $\{e_j\}$ are orthonormal vectors in $\HH$ and $j_{\ell} \ne j_{\ell+1}$ for $1 \le \ell \le k-1$, then
\begin{equation}\label{action of Sq}
{\mathscr S}_q\colon H^q_{n_1}(X_q(e_{j_1}))\cdots H^q_{n_k}(X_q(e_{j_k}))\mapsto Z_q(e_{j_1})^{n_1}\cdots Z_q(e_{j_k})^{n_k}.
\end{equation}
Comparing Equations \ref{action of S} and \ref{action of Sq}, we see that ${\mathscr S}_q$ is a natural extension of the classical Segal-Bargmann transform.

\medskip

Since $\H_q(\Hc)$ is contained in $\Gamm_q(\HH\oplus\HH)$, it inherits the number operator $N_q$ from it, induced by the inclusion of $\Gamm_q(\HH\oplus\HH)$ into $\F_q(\HH\oplus\HH)$ via the map $A\mapsto A\vac$.  From Equation \ref{monomial mapping}, we see then that $Z_q(e_1)^{n_1}\cdots Z_q(e_k)^{n_k}$ is an eigenvector of $N_q$ with eigenvalue $n_1+\cdots+n_k$.  This is precisely matches the conjugated action ${\mathscr S}_qN_q{\mathscr S}_q^{\ast}$ of the number operator $N_q$ on $\Gamm_q(\HH)$, as can be seen from Proposition \ref{Wick ordering 1}. Hence, we have ${\mathscr S}_qN_q = N_q{\mathscr S}_q$, just as in the commutative case.

\medskip

Finally, we define $L^p(\H_q(\Hc),\tau_q)$ to be the completion of $\H_q(\Hc)$ in the $L^p(\Gamm_q(\HH\oplus\HH),\tau_q)$-norm. For $p\ge 2$ (the case of interest for our main theorem), it is equal to the intersection of $L^2(\H_q(\Hc),\tau_q)$ with $L^p(\Gamm_q(\HH\oplus\HH),\tau_q)$.  The class of Banach spaces $L^p(\H_q,\tau_q)$ is a non-commutative generalization of the spaces $\H L^p(\C^d,\gamma')$ that occur in Janson's Theorem \ref{Janson hyp thm}.  Since the algebra $\H_q$ is not a von Neumann algebra, this family is {\it not known} to be complex interpolation scale.  For example, in the $q=1$ case, the family is not complex interpolation scale when $\HH$ is infinite-dimensional (this is almost proven in \cite{JPR}).  Hence, once we have proved Theorem \ref{Main Theorem}, it is not an easy matter to generalize to the case $p> 2$, $r\ne 2,4,6,\ldots$

\section{Mixed Spin and Strong Hypercontractivity}\label{spin alg}

We will consider the mixed-spin algebras $\CC(I,\eps)$ introduced in \cite{Biane 1} which represent systems with some commutation and some anti-commutation relations.  Such systems may be viewed as approximations to the $q$-commutation relations, in a manner which will be made precise
in Section \ref{Speicher}.  We introduce a holomorphic subalgebra $\H(I,\eps)$, and give a combinatorial proof of a strong hypercontractivity theorem like Theorem \ref{Main Theorem} for it.

\subsection{The mixed-spin algebra $\CC(I,\eps)$}  Let $I$ be a finite totally ordered set (with cardinality denoted by $|I|$), and let $\eps$ be a function $I\times I\to\{-1,1\}$ which is symmetric, $\eps(i,j) = \eps(j,i)$, and constantly $-1$ on the diagonal, $\eps(i,i) = -1$.  Let $\CC(I,\eps)$ denote the unital $\C$-algebra with generators $\{x_i\,;\,i\in I\}$ and relations
\begin{equation}\label{eps relations}
x_i x_j - \eps(i,j) x_j x_i = 2\delta_{ij} \quad\text{for}\quad i,j\in I.
\end{equation}
(The requirement $\eps(i,i) = -1$ forces $x_i^2 = 1$, and guarantees that $\CC(I,\eps)$ is finite-dimensional.) In the special case $\eps\equiv -1$, this is precisely the complex Clifford algebra $\CC_{|I|}$, hence our choice of notation.  In the case $\eps(i,j) = 1$ for $i\ne j$ (i.e.\ when different generators commute), the generators of $\CC(I,\eps)$ may be modeled by $|I|$ i.i.d.\ Bernoulli random variables, and so we reproduce the toy Fock space considered in \cite{Meyer}.  In the general case, $\CC(I,\eps)$ has, as a vector space, a basis consisting of all $x_A$ with $A = (i_1,\ldots,i_k)$ increasing multi-indices in $I^k$, where $x_A = x_{i_1}\ldots x_{i_k}$, and $x_{\emptyset}$ denotes the identity $1\in\CC(I,\eps)$.  Thus, $\dim\CC(I,\eps) = 2^{|I|}$.  Moreover, $\CC(I,\eps)$ has a natural
decomposition
\[ \CC(I,\eps) = \bigoplus_{n=0}^{|I|} \CC_n(I,\eps), \]
where $\CC_n = \text{span}\{x_A\,;\,|A| = n\}$ is the ``$n$-particle space."  Of some importance to us will be the natural grading of
the algebra,
\[ \CC(I,\eps) = \CC_+(I,\eps)\oplus\CC_-(I,\eps), \]
where $\CC_+ = \bigoplus\{\CC_n\,;\,n\text{ is even}\}$, and $\CC_-$ is the corresponding odd subspace.  The reader may readily verify
that this decomposition is a grading -- i.e.\ $\CC_\alpha\cdot\CC_\beta\subseteq\CC_{\alpha\beta}$, where $\alpha,\beta\in\{+,-\}$ and their
product is to be interpretted in the obvious fashion.

\medskip

We equip $\CC(I,\eps)$ with an involution $\ast$, which is defined to be the conjugate-linear extension of the map $x_A^\ast = x_{A^\ast}$,
where $(i_1,\ldots,i_k)^\ast$ is the reversed multi-index $(i_k,\ldots,i_1)$.  In particular, the generators $x_i = x_i^\ast$ are self-adjoint,
and in general $x_A^\ast = \pm x_A$.  We also define a tracial state $\te$ by $\te(x_A) = \delta_{A\emptyset}$; that is, $\te(1) = 1$ while $\te(x_A) = 0$
for all other basis elements.  It is easy to check that $\te(ab) = \te(ba)$.  This allows us to define an inner product on $\CC(I,\eps)$ by
\[ (a,b)_\eps = \te(b^\ast a). \]
The basis $\{x_A\}$ is orthonormal with respect to $(\cdot,\cdot)_\eps$.  Following the GNS construction, the action of $\CC(I,\eps)$ on the Hilbert space $(\CC(I,\eps),(\cdot,\cdot)_\eps)$ by left-multiplication is continuous, and yields an injection of $\CC(I,\eps)$ into the von Neumann algebra
of bounded operators on the Hilbert space.  In this way, $\CC(I,\eps)$ gains a von Neumann algebra structure.  We denote by $L^p(\CC(I,\eps),\te)$
the non-commutative $L^p$ space of this von Neumann algebra with its trace $\te$.  (So, in particular, $L^2(\CC(I,\eps),\te)$ is naturally isomorphic to the Hilbert space $(\CC(I,\eps),(\cdot,\cdot)_\eps)$.) The $L^p(\CC(I,\eps),\te)$-norm is, in fact, just the (normalized) Schatten $L^p$-norm on the matrix algebra.  This can be seen from the following Proposition.
\begin{proposition}\label{trace}
Let $tr$ denote the normalized trace on the finite-dimensional algebra $\B(L^2(\CC(I,\eps),\te))$.  Then for any $x\in\CC(I,\eps)$,
\[ \te(x) = tr(x). \]
\end{proposition}
\begin{proof}
Using the orthonormal basis $\{x_A\}$ for $L^2(\CC(I,\eps),\te)$, we compute that
\[ tr(x_A) = 2^{-|I|}\sum_B (x_Ax_B,x_B)_\eps = \begin{cases} 0,&\text{ if }A\ne\emptyset \\ 2^{-|I|}\sum_B(x_B,x_B)_\eps = 1,&\text{ if }A=\emptyset \end{cases} \]
where the sums are taken over all increasing multi-indices $B$.  So $tr(x_A) = \delta_{A\emptyset} = \te(x_A)$.  
\end{proof}
Note that the trace $\te$ may be expressed in terms of the inner product as the pure state $\te(x) = (x1,1)_\eps$.  This formula extends to all of $\B(L^2(\CC(I,\eps),\te))$, giving the pure state $\beta\mapsto(\beta 1,1)_\eps$.  However, this state does {\it not} equal $tr$ for all bounded operators $\beta$.  We will see examples in Section \ref{Speicher} showing that it is not tracial in general.

\medskip

\noindent The algebra $\CC(I,\eps)$ comes equipped with a number operator $N_\eps$ which has $\CC_n(I,\eps)$ as an eigenspace with eigenvalue $n$.  That is, \[ N_\eps x_A = |A|x_A. \] This is a generalization of the action of the operator $N_{-1}$ on the Clifford algebra $\CC_{|I|}=\Gamm_{-1}(\R^{|I|})$. $N_\eps$ is a positive semi-definite operator on $L^2(\CC(I,\eps),\te)$, and so generates a contraction semigroup $e^{-tN_\eps}$.  It is to the study of this semigroup, restricted to a holomorphic subspace, that we devote the remainder of this section.

\subsection{The mixed-spin holomorphic algebra $\H(I,\eps)$}\label{holomorphic algebra}

Following our construction of $\H_q$, we will begin by doubling the number of variables.  We extend $\eps$ to the set $I\times\{0,1\}$ by setting
\[ \eps((i,\zeta),(j,\zeta')) = \eps(i,j), \]
and then consider the algebra $\CC(I\times\{0,1\},\eps)$.  If we relabel $x_{(i,0)}\to x_i$ and $x_{(i,1)}\to y_i$, then this is tantamount to constructing
the unital $\C$-algebra with relations
\begin{equation}\label{eps relations 2} \left.\begin{matrix} x_i x_j - \eps(i,j) x_j x_i & = & 2\delta_{ij} \\
                             y_i y_j - \eps(i,j) y_j y_i & = & 2\delta_{ij} \\
                             x_i y_j - \eps(i,j) y_j x_i & = & 0 \end{matrix}\right\}\quad \text{for} \quad i,j\in I.
\end{equation}
Note, $\CC(I,\eps)$ is $\ast$-isomorphically embedded in $\CC(I\times\{0,1\},\eps)$ via the inclusion $x_i\mapsto x_{(i,0)}$.  Hence, this relabeling should not be confusing.

\medskip

We define elements $z_j\in\CC(I\times\{0,1\},\eps)$ by
\begin{equation}\label{z eps}
z_j = 2^{-1/2}(x_j + iy_j) = 2^{-1/2}(x_{(j,0)} + ix_{(j,1)}),\quad \text{for} \quad j\in I.
\end{equation}
(To avoid confusion, we point out that in Equation \ref{z eps}, $i$ refers to $\sqrt{-1}\in\C$.)  The operator $z_j$ is an analog of the operators
$Z_q(e_j)$ in $\H_q$.  The normalization is again chosen so that $z_j$ is a unit vector in $L^2(\CC,\te)$.  For the calculations in the foregoing, however,
it will be convenient to have the variables normalized in $L^\infty(\CC,\te)$.  Therefore, we also introduce
\[ \z_j = 2^{-1/2} z_j = \frac{1}{2}(x_j + i y_j)\quad \text{for} \quad j\in I. \]
The reader may readily verify that $|\z_j|^2 = \z_j^\ast\z_j$ is a nonzero idempotent, and hence $\|\z_j\|_\infty = 1$.  Define the {\bf mixed spin holomorphic algebra} $\H(I,\eps)$ as the $\C$-algebra generated by $\{z_1,\ldots,z_{|I|}\}$.  This is just the polynomial algebra in the variables $z_1,\ldots,z_{|I|}$ --- the adjoints are not included.  Indeed, $2\z_j^\ast = x_j - iy_j$, so $x_j = \z_j + \z_j^\ast$ and $y_j = i(\z_j^\ast - \z_j)$.  Thus, the $\ast$-algebra
generated by $z_1,\ldots,z_{|I|}$ is all of $\CC(I\times\{0,1\},\eps)$.

\medskip

Observe that $2z_j^2 = x_j^2 - y_j^2 + i(x_j y_j + y_j x_j) = 0$ since $\eps(j,j) = -1$.  In general, we may compute that
\begin{equation}\label{z commute}
z_i z_j - \eps(i,j) z_j z_i = 0,
\end{equation}
and the same relations (of course) hold for the $\z_j$.  The operators $\z_j,\z_j^\ast$  also satisfy the joint relations
\begin{equation}\label{eps ast relations}
\z_i^\ast \z_j - \eps(i,j) \z_j \z_i^\ast = \delta_{ij} \quad \text{for} \quad i,j\in I.
\end{equation}
Equation \ref{eps ast relations} looks much like the $q$-commutation relations of Equation \ref{q-commutation relations}.  It is, in fact, possible to think of $\z_j,\z_j^\ast$ as creation and annihilation operators.  That is, there is a faithful representation of $\z_j,\z_j^\ast$ {\it in} $\B(L^2(\CC(I,\eps),\te))$, which sends $\z_j$ and $\z_j^\ast$ to the creation and annihilation operators $\beta_j,\beta_j^\ast$ on $L^2(\CC(I,\eps),\te)$ discussed in Section \ref{Speicher}.  (By our definition, the operators $\z_j,\z_j^\ast$ are a priori in the {\it doubled} space $\B(L^2(\CC(I\times\{0,1\},\eps),\te)$.)  This representation, the {\it spin-chain representation}, is discussed in \cite{Carlen Lieb} in detail in the case $\eps\equiv -1$, and is generalized in \cite{Biane 1}.  The problem with this point of view is that the pure state $\beta\mapsto(\beta1,1)_\eps$ on $\B(L^2(\CC(I,\eps),\te))$ does {\it not} correspond to the trace $\te$ under the representation.  So, we prefer not to think of $\z_j,\z_j^\ast$ as creation and annihilation operators.

\medskip

A simple calculation shows that if $|A| = n$ then $z_A \in \CC_n(I\times\{0,1\},\eps)$, and so $N_\eps z_A = |A|z_A$.  Thus, $\H(I,\eps)$ is a reducing subspace for the (self-adjoint) operator $N_\eps$ on $L^2(\Cd,\te)$.  Note also that the action of $N_\eps$ on $z_A$ mirrors that of $N_\eps$ on $x_A$.
In fact, this can be stated in terms of a {\bf $\eps$-Segal-Bargmann transform}: the map ${\mathscr S}_\eps\colon x_A\mapsto z_A$
is a unitary isomorphism of $L^2(\Cs,\te)$ onto $L^2(\H(I,\eps),\te)$, and ${\mathscr S}_\eps N_\eps = N_\eps{\mathscr S}_\eps$.

\medskip

The main part of the proof of Theorem \ref{Main Theorem} is the following strong hypercontractivity result regarding the semigroup
$e^{-tN_\eps}$ acting on $\H(I,\eps)$.
\begin{theorem}\label{eps strong hyp} For $p=2$ and $r$ an even integer,
\[ \|e^{-tN_\eps} a\|_r \le \|a\|_p\text{ for all }a\in\H(I,\eps) \quad \text{iff} \quad t\ge t_J(p,r) = \frac{1}{2}\log\frac{r}{p}. \] 
\end{theorem}
We expect the theorem holds for $2\le p \le r < \infty$. (The case $p<2$ may be somewhat different from the commutative case; in a communication from L.\ Gross, a calculation showed that in $1$-dimension the least time to contraction seems to be larger than the Janson time for some $p,r<2$.)  If $I\ne\emptyset$, it is easy to see that the Janson time cannot be improved for {\it any $p,r>0$}, again by calculation in the $1$-dimensional case.
\begin{proof}[Proof of the `only if' direction of theorem \ref{eps strong hyp}] Let $a(\e) = 1+\e \z\in \H(I,\eps)$ where $\z = \z_j$ for some
$j\in I$.  Then $|a(\e)|^2 = (1+\e\z^\ast)(1+\e\z) = 1 + \e(\z + \z^\ast) + \e^2|\z|^2 = 1 + \e x + \e^2|\z|^2$, where $x = x_j$.
Hence,
\begin{eqnarray*}
|a(\e)|^{2p} & = & (1+\e(x + \e|\z|^2))^p \\
                      & = & 1 + p\e(x + \e|\z|^2) + \frac{p(p-1)}{2}\e^2(x+\e|\z|^2)^2 + o(\e^2) \\
                      & = & 1 + \e(px) + \e^2\left(p|\z|^2 + \frac{p(p-1)}{2}x^2\right) + o(\e^2).
\end{eqnarray*}
Now, $|\z|^2 = (1/2)(1+ixy)$ where $y=y_j$, and so $\te|\z|^2 = 1/2$.  Also $x^2 = 1$, and $\te x = 0$.  Therefore,
\begin{eqnarray*}
\|a(\e)\|_{2p} \;\; = \;\; \left(\te(|a(\e)|^{2p})\right)^{1/2p} & = & \left(1 + \left(p\cdot\frac{1}{2} + \frac{p(p-1)}{2}\right)\e^2 + o(\e^2)\right)^{1/2p} \\
& = & 1 + \frac{p}{4}\e^2 + o(\e^2).
\end{eqnarray*}
So $\|a(\e)\|_p = 1 + (p/8)\e^2 + o(\e^2)$.  Now, $e^{-tN_\eps}a(\e) = 1 + \e e^{-t}\z = a(e^{-t}\e)$.  Thus, in order
for $\|e^{-tN_\eps}a(\e)\|_r \le \|a(\e)\|_p$, we must have
\[ 1+ \frac{1}{8}e^{-2t}r\e^2 + o(\e^2) \le 1 + \frac{1}{8}p\e^2 + o(\e^2), \]
and so as $\e\to 0$, it follows that $e^{-2t}\le p/r$ --- or $t \ge t_J(p,r)$.  
\end{proof}
Hence, the necessity condition holds for all $r\ge p > 0$.  For the sufficiency, however, the tools available to us are extremely limited (due to the
fact that $\H(I,\eps)$ is not a $\ast$-algebra).  We are forced to give a combinatorial proof, which cannot reach beyond the cases when $p=2$ and $r$ is even.  The remainder of the `if' direction of Theorem \ref{eps strong hyp} is the main subject of Section \ref{H strong hyp}.

\subsection{Strong hypercontractivity for $\H(I,\eps)$}\label{H strong hyp}
We will prove Theorem \ref{eps strong hyp} by induction on $|I|$.  Note, in the case $|I|=0$, the algebra $\H(I,\eps)$ is just $\C$.  Since the action of $e^{-tN_\eps}$ on $\C$ is trivial, and since all $\|\cdot\|_p$ norms are equal to the complex modulus $|\cdot|$, the sufficiency condition follows automatically in this case. Now, suppose the strong hypercontractivity result of Theorem \ref{eps strong hyp} holds for the algebras $\H(I',\eps')$ with $|I'|\le d$.  Let $I$ be a set of size $d+1$, and $\eps$ a spin-assignment on $I$.  Select any fixed element $i\in I$.  Any element $a\in\H(I,\eps)$ can be uniquely decomposed as
\begin{equation}\label{decomp}
a = b + \z_ic, \quad b,c\in\H(I-\{i\},\left.\eps\right|_{I-\{i\}}).
\end{equation}
For convenience, throughout we will refer to $I-\{i\}$ as $J$, and to $\z_i$ as $\z$, $x_i$ as $x$, and so forth.  Since $|J| = d$, the inductive hypothesis
is that $\H(J,\left.\eps\right|_J)$ satisfies the strong hypercontractivity estimate of Theorem \ref{eps strong hyp}.

\medskip

The quantity $|\z|^2$ will often come up in calculations, and so we give it a name: $\xi = |\z|^2 = \z^\ast\z$.  We will also encounter $\z\z^\ast$,
but by Equation \ref{eps ast relations}, $\z\z^\ast = 1-\xi$.  The following lemma records some of the important properties of the operators $\xi$, $\z$, and $\z^\ast$.  All of the statements may be verified by trivial calculation.
\begin{lemma}\label{calc}  The following properties hold for $\xi$, $\z$, and $\z^\ast$.
\begin{enumerate}
\item $\xi^p = \xi$ for $p > 0$.
\item $\xi$ is independent of $\CC(J\times\{0,1\},\left.\eps\right|_J)$ --- that is, for each $u\in\CC(J\times\{0,1\},\left.\eps\right|_J)$, $\xi u = u\xi$ and $\te(\xi u) = \te(\xi)\te(u) = \frac{1}{2}\te(u)$.\label{independence}
\item Let $u\in\CC(J\times\{0,1\},\left.\eps\right|_J)$, let $h\in\{\z,\z^\ast,\xi,1-\xi\}$, and let $p>0$.  Then $\|hu\|_p = 2^{-1/p}\|u\|_p$.\label{product norms}
\item $\xi\z = \z^\ast\xi = 0$, $\z\xi = \z$, and $\xi\z^\ast = \z^\ast$\label{products}.
\end{enumerate}
\end{lemma}
The commutativity in item \ref{independence} above follows in large part from the fact that $\xi = \frac{1}{2}(1 + ixy) \in \CC_+(I\times\{0,1\},\eps)$.
The grading plays an important role in the combinatorics to follow.  In fact, the grading of $\CC(\{i\}\times\{0,1\},\left.\eps\right|_{\{i\}})$ induces a grading
on the full algebra $\CC(I\times\{0,1\},\eps)$.  We refer to this grading by
\[ \CC = \CC_+^i\oplus\CC_-^i, \qquad \CC_\alpha^i\cdot\CC_\beta^i \subseteq \CC_{\alpha\beta}^i. \]
So, for example, the element $y_j\xi$ ($i\ne j$) is in $\CC_+^i(I\times\{0,1\},\eps)$, even though it is in $\CC_-(I\times\{0,1\},\eps)$.  Note that
\[ \CC^i_-(I\times\{0,1\},\eps) = \left\{\z u + \z^\ast v\,;\, u,v\in \CC(J\times\{0,1\},\left.\eps\right|_J)\right\}. \]  For any such $u$, $\te(\z u) = (\z,u^\ast)_\eps = 0$,
and similarly $\te(\z^\ast u) = 0$.  It follows that $\te|_{\CC_-^i} = 0$.  Using the graded structure, this leads to the following important lemma, which
aids in the calculation of moments.
\begin{lemma}\label{even}
Let $v^0\in \CC^i_+(I\times\{0,1\},\eps)$ and $v^1\in \CC^i_-(I\times\{0,1\},\eps)$.  Let $\eta$ be $\{0,1\}$-sequence of length $n$, and denote by $|\eta|$ the sum $\eta_1+\cdots+\eta_n$ of its entries (i.e.\ the number of $1$s). Then the element $v^\eta = v^{\eta_1}\cdots v^{\eta_n}$ has $\te(v^\eta) = 0$ if $|\eta|$ is odd.
\end{lemma}
Now, we proceed to expand the moments of $|a|^2$.  Using the decomposition in \ref{decomp}, we have $|a|^2 = (b+\z c)^\ast(b+\z c) = |b|^2 + b^\ast \z c + c^\ast\z^\ast b + c^\ast|\z|^2c$.  That is,
\begin{equation}\label{a^2}
|a|^2 = (|b|^2 + \xi |c|^2) + (b^\ast \z c + c^\ast \z^\ast b) = v^0 + v^1.
\end{equation}
Equation \ref{a^2} decomposes $|a|^2$ into its $\CC^i_+$ and $\CC^i_-$ parts, $v^0 = |b|^2 + \xi |c|^2$ and $v^1 = b^\ast \z c + c^\ast \z^\ast b$.  It follows immediately that
\begin{equation}\label{a 2 norm}
\|a\|_2^2 = \te(v^0) = \te(|b|^2 + \xi |c|^2) = \|b\|_2^2 + \frac{1}{2}\|c\|_2^2.
\end{equation}
The factor of $1/2$ (unusual in Pythagoras' formula) is due to our choice to normalize $\z$ in $L^\infty$ and not in $L^2$.  More generally, for
the $n$th moment of $|a|^2$,
\[ \|a\|_{2n}^{2n} = \te(|a|^{2n}) = \te[(v^0 + v^1)^n] = \sum_{\eta\in 2^n} \te(v^\eta), \]
where $2^n$ denotes the set of all $\{0,1\}$-sequences of length $n$.  Using Lemma \ref{even}, we have
\[ \|a\|_{2n}^{2n} = \sum_{\stackrel{|\eta|\text{ even}}{\eta\in 2^n}} \te(v^\eta) = \sum_{k=0}^{\lfloor n/2 \rfloor} \sum_{|\eta|=2k} \te(v^\eta). \]
Now, the term $v^\eta$ is a product of $n$ terms, each of which is either $|b|^2 + \xi |c|^2$ or $b^\ast \z c + c^\ast \z^\ast b$.  Define
\[ \begin{aligned}
v^{00} &= |b|^2 \qquad &v^{01} &= \xi |c|^2 \\
v^{10} &= b^\ast \z c \qquad &v^{11} &= c^\ast \z^\ast b.
\end{aligned} \]
Then we may write $v^\eta$ as
\[ v^\eta = \sum_{\nu\in 2^n} v^{\eta\nu} =  \sum_{\nu\in 2^n} v^{\eta_1\nu_1}\cdots v^{\eta_n\nu_n}. \]

\medskip

It should be noted that many of the terms in this sum are in fact $0$.  For example, consider $(v^{10})^2 = b^\ast \z c b^\ast \z c$.  In general, for any $u\in\CC(J\times\{0,1\},\left.\eps\right|_J)$, there is a $\tilde{u}\in\CC(J\times\{0,1\},\left.\eps\right|_J)$ such that $\z u = \tilde{u}\z$.  Hence the term $(v^{10})^2$ contains $\z^2 = 0$, and so is $0$.  More generally,
a term like $v^{10} v^{01} v^{10}$ is also $0$: the $\z$ in $v^{10}$ can be commuted past all terms except $\xi$, at which point the product is either
$0$ or $\z$ (by Lemma \ref{calc}), so the term is $0$.  On the other hand, the term $v^{11}v^{01}v^{10}$ is nonzero, since (once commuting past the
$\CC(J\times\{0,1\},\left.\eps\right|_J)$-terms) we have $\z^\ast\xi\z = (\z^\ast \z)^2 = \xi \ne 0$.

\medskip

Let $\eta,\nu\in 2^n$.  Denote by $\1(\eta)\subseteq\{1,\ldots,n\}$ the set of $j$ such that $\eta_j = 1$.  Then say that $\nu$ is {\bf $\eta$-alternating},
$\nu\in A(\eta)$, if the subsequence $\{(\nu_j)\,;\,j\in\1(\eta)\}$ is alternating.  For example, let $\eta = (1,1,0,1)$.  Then the sequences
$(0,1,0,0)$ and $(0,1,1,0)$ are both in $A(\eta)$, while the sequence $(0,0,0,0)$ is not.  Note that $v^{10}$ and $v^{11}$ are the terms containing
$\z$ and $\z^\ast$.  Hence, the $v^{\eta\nu}$ with $\nu\in A(\eta)$ are precisely those terms in which $\z$ and $\z^\ast$ alternate when they occur.
By the considerations in the preceeding paragraph, these are the only nonzero terms in the expansion of $v^\eta$. Thus,
\[ v^\eta = \sum_{\nu\in A(\eta)} v^{\eta\nu}. \]
In any term in the above sum, let $|\eta|=2k$ and let $|\nu|=m$.  Since $\1(\eta)$ is a set of $2k$ indices, and since $\nu\in A(\eta)$, $\nu_j=0$ for $k$
of these indices $j$, and $\nu_j=1$ for the other $k$.  Thus $\nu$ contains at least $k$ $1$s and at least $k$ $0$s, and so $k \le m \le n - k$.  It follows that the full expansion for the $n$th moment is
\[ \|a\|_{2n}^{2n} = \sum_{k=0}^{\lfloor n/2 \rfloor}\sum_{m=k}^{n-k} \sum_{|\eta|=2k} \sum_{\stackrel{\nu\in A(\eta)}{|\nu|=m}} \te(v^{\eta\nu}). \]
It will be useful to consider the cases $k=0$ and $m=0$ separately, and so we rewrite this moment as
\begin{equation}\label{moment}
\|a\|_{2n}^{2n} = \te[(v^{00})^n] + \sum_{m=1}^n \sum_{|\nu|=m} \te(v^{0\nu}) + \sum_{k=1}^{\lfloor n/2 \rfloor}\sum_{m=k}^{n-k} \sum_{|\eta|=2k} \sum_{\stackrel{\nu\in A(\eta)}{|\nu|=m}} \te(v^{\eta\nu}).
\end{equation}
(Note, if $\eta\equiv 0$ then the condition $\nu\in A(\eta)$ is vacuously satisfied for all $\nu\in 2^n$.)  Each of the $v^{\eta\nu}$ in Equation
\ref{moment} is a product of terms, each of which contains some elements of $\CC(J\times\{0,1\},\left.\eps\right|_J)$ and some factors of $\z$, $\z^\ast$, or $\xi$. (Observe the only term which has no factors from $\CC(\{i\}\times\{0,1\},\left.\eps\right|_{\{i\}})$ is  the first one $(v^{00})^n$.)  To estimate such terms, we introduce the following tool.
\begin{lemma}\label{estimate}
Let $u_1,\ldots,u_s \in \CC(J\times\{0,1\},\left.\eps\right|_J)$.  Let $U$ be a product including all of the elements $u_1,\ldots,u_n$ together with some non-zero number of terms from $\{\z,\z^\ast,\xi\}$.  Then
\begin{equation}\label{estimate eqn}
\te(U) \le \frac{1}{2}\|u_1\|_s\cdots\|u_s\|_s.
\end{equation}
\end{lemma}
\begin{proof} First note that $\te(U)$ is invariant under cyclic permutations of $U$.  $U$ may then be written in the form $h_1U_1h_2U_2\cdots h_\ell U_\ell$, where each $U_j$ is a product of some of the $u_1,\ldots,u_s$, and each $h_j$ is a product of the terms $\z$,$\z^\ast$, and $\xi$.  Let $s_j$ be the number of terms in $U_j$; then $s_1+\cdots + s_\ell = s$.  So $s_1/s + \cdots + s_\ell/s = 1$, and when we apply H\"older's inequality, we find
\begin{equation}\label{Holder 1}
\te(U) \le \|h_1 U_1\|_{s/s_1}\cdots \|h_\ell U_\ell\|_{s/s_\ell}.
\end{equation}
By Lemma \ref{calc} part \ref{products}, any product of terms in $\{\z,\z^\ast,\xi\}$ is either $\z$, $\z^\ast$, $\xi$, $1-\xi$, or $0$.  Thus, using Lemma \ref{calc} part \ref{product norms}, we have
\begin{equation}\label{Holder 2}
\|h_j U_j\|_{s/s_j} \le 2^{-s_j/s}\|U_j\|_{s/s_j}.
\end{equation}
Now, since $U_j$ is a product of $s_j$ terms, say $u_{k_1},\ldots,u_{k_{s_j}}$, applying H\"older's inequality again (using $1/s_j + \cdots + 1/s_j = 1/(s/s_j)$) we have $\|U_j\|_{s/s_j} \le \|u_{k_1}\|_s\cdots\|u_{k_{s_j}}\|_s$.  Combining this with Equations \ref{Holder 1} and \ref{Holder 2}, we get
\[ \te(U) \le 2^{-s_1/s}\cdots 2^{-s_\ell/s} \|u_1\|_s \cdots \|u_s\|_s, \]
and since $s_1+\cdots+s_\ell = s$, this reduces to Equation \ref{estimate eqn}.  
\end{proof}
We now apply Lemma \ref{estimate} to estimate the three terms in Equation \ref{moment}.  The first term is merely $\te(|b|^{2n}) = \|b\|_{2n}^{2n}$.
In the first sum
\[ \sum_{m=1}^n \sum_{|\nu|=m} \te(v^{0\nu}), \]
the term $v^{0\nu}$, with $|\nu| = m$, is a product containing $m$ factors of $v^{01} = \xi c^\ast c$ and $n-m$ factors of $v^{00} = b^\ast b$.
So there are a total of $2n$ factors from the set $\{b,b^\ast,c,c^\ast\}\subset\CC(J\times\{0,1\},\left.\eps\right|_J)$.  since $\|u\|_{2n} = \|u^\ast\|_{2n}$ for each $u\in\CC(J\times\{0,1\},\left.\eps\right|_J)$, Lemma \ref{estimate} then implies that
\[ \te(v^{0\nu}) \le \frac{1}{2}(\|b\|_{2n})^{2(n-m)}(\|c\|_{2n})^{2m}. \]
Hence,
\begin{eqnarray}
\sum_{m=1}^n \sum_{|\nu|=m} \te(v^{0\nu}) & \le & \frac{1}{2}\sum_{m=1}^n \sum_{|\nu|=m}(\|b\|_{2n}^2)^{n-m}(\|c\|_{2n}^2)^m \nonumber\\
& = & \frac{1}{2}\sum_{m=1}^n\binom{n}{m} (\|b\|_{2n}^2)^{n-m}(\|c\|_{2n}^2)^m. \label{sum 1}
\end{eqnarray}
Now we consider the second sum
\[ \sum_{k=1}^{\lfloor n/2 \rfloor}\sum_{m=k}^{n-k} \sum_{|\eta|=2k} \sum_{\stackrel{\nu\in A(\eta)}{|\nu|=m}} \te(v^{\eta\nu}). \]
In each term $v^{\eta\nu}$, since $|\eta|=2k$ and $\nu\in A(\eta)$, we know that $k$ of the terms are $v^{10}$ and $k$ of the terms are $v^{11}$.  So $k$ of the $1$s in $\nu$ have been accounted for with the $v^{11}$ terms, and since $|\nu|=m$ precisely $m-k$ terms must be $v^{01}$. As the total number of terms must be $n$, this means the remaining $v^{00}$ terms are $n-(2k+m-k) = n-m-k$ in number.  So, there are
\begin{itemize}
\item $k$ factors of $v^{10} = b^\ast \z c$, so $k$ factors each of $b^\ast$ and $c$,
\item $k$ factors of $v^{11}= c^\ast \z b$, so $k$ factors each of $b$ and $c^\ast$,
\item $m-k$ factors of $v^{01} = \xi c^\ast c$, so $m-k$ factors each of $c$ and $c^\ast$, and
\item $n-m-k$ factors of $v^{00} = b^\ast b$, so $n-m-k$ factors each of $b$ and $b^\ast$.
\end{itemize}
In total, then, $v^{\eta\nu}$ contains $2k + 2(n-m-k) = 2(n-m)$ factors of $b$ or $b^\ast$, and $2k + 2(m-k) = 2m$ factors of $c$ or $c^\ast$.  Applying Lemma \ref{estimate} again,
\[ \sum_{k=1}^{\lfloor n/2 \rfloor}\sum_{m=k}^{n-k} \sum_{|\eta|=2k} \sum_{\stackrel{\nu\in A(\eta)}{|\nu|=m}} \te(v^{\eta\nu})
\le \frac{1}{2}\sum_{k=1}^{\lfloor n/2 \rfloor}\sum_{m=k}^{n-k} \sum_{|\eta|=2k} \sum_{\stackrel{\nu\in A(\eta)}{|\nu|=m}} (\|b\|_{2n})^{2(n-m)}(\|c\|_{2n})^{2m}.
\]
We must now count the number of pairs $(\eta,\nu)$ with $|\eta|=2k$, $\nu\in A(\eta)$ and $|\nu|=m$.  There are $\binom{n}{2k}$ such $\eta$.  We know that $\nu$ is alternating on $\1(\eta)$, and so the corresponding subsequence must be either $0101\ldots01$ or $1010\ldots10$, giving two choices, and exhausting $k$ of the $m$ $1$s in $\nu$.  Finally, since $|\1(\eta)| = 2k$, there are $n-2k$ $0$s in $\eta$, and $\nu$ is unconstrained there; hence, there are $\binom{n-2k}{m-k}$ choices. Whence, the number of pairs $(\eta,\nu)$ in the sum is
\[ 2\binom{n}{2k}\binom{n-2k}{m-k}. \]
This gives the estimate
\begin{eqnarray}
&& \sum_{k=1}^{\lfloor n/2 \rfloor}\sum_{m=k}^{n-k} \sum_{|\eta|=2k} \sum_{\stackrel{\nu\in A(\eta)}{|\nu|=m}} \te(v^{\eta\nu}) \nonumber \\
& \le & \sum_{k=1}^{\lfloor n/2 \rfloor}\sum_{m=k}^{n-k} \binom{n}{2k}\binom{n-2k}{m-k} (\|b\|_{2n}^2)^{n-m}(\|c\|_{2n}^2)^m \label{sum 2}
\end{eqnarray}
for the final sum in Equation \ref{moment}.  It will be convenient to reorder the terms in Equation \ref{sum 2} so that $m$ occurs first.  Since the sum (for each $k$) ranges from $m=k$ to $m=n-k$, the pairs $(k,m)$ in the sum are those with $1\le k\le\lfloor n/2\rfloor$ and $k\le m\le n-k$.  The second condition gives two inequalities: $k\le m$ and $k\le n-m$.  Note, if both of these are satisfied then, summing, $2k\le n$ -- the first condition is automatically satisfied.  The
sum can therefore be rewritten as
\begin{equation}\label{sum 3}
\sum_{m=1}^n\sum_{k=1}^{m\wedge(n-m)}\binom{n}{2k}\binom{n-2k}{m-k}(\|b\|_{2n}^2)^{n-m}(\|c\|_{2n}^2)^m.
\end{equation}
So, combining Equations \ref{moment}, \ref{sum 1}, and \ref{sum 3}, we have the estimate
\begin{equation}\label{moment estimate}
\|a\|_{2n}^{2n} \le \|b\|_{2n}^{2n} + \sum_{m=1}^n \chi_m (\|b\|_{2n}^2)^{n-m}(\|c\|_{2n}^2)^m,
\end{equation}
where the coefficient $\chi_m$ is given by
\[ \chi_m = \frac{1}{2}\binom{n}{m} + \sum_{k=1}^{m\wedge(n-m)}\binom{n}{2k}\binom{n-2k}{m-k}. \]
The following proposition shows that $\chi_m$ is optimally bounded to yield the necessary strong hypercontractive estimate.  We state it without proof;
the reader may do the necessary calculations.
\begin{proposition}\label{combinatorics}
The coefficients $\chi_m$ satisfy \[ \chi_m \le \binom{n}{m}\left(\frac{n}{2}\right)^m. \]  This inequality is an equality in the case $m=1$.
\end{proposition}
Applying Proposition \ref{combinatorics} to Equation \ref{moment estimate}, we have
\begin{equation}\label{moment estimate 2}
\|a\|_{2n}^{2n} \le \|b\|_{2n}^{2n} + \sum_{m=1}^n \binom{n}{m}\left(\frac{n}{2}\right)^m (\|b\|_{2n}^2)^{n-m}(\|c\|_{2n}^2)^m.
\end{equation}
We now complete the proof of Theorem \ref{eps strong hyp}.
\begin{proof}[Proof of the `if' direction of Theorem \ref{eps strong hyp}] For $a = b + \z c$, we have $a_t = b_t + e^{-t}\z c_t$, where $a_t = e^{-tN_\eps}a$ and so forth.  Using the estimate in Equation \ref{moment estimate 2}, we have
\[ \|a_t\|_{2n}^{2n} \le \|b_t\|_{2n}^{2n} + \sum_{m=1}^n \binom{n}{m}\left(\frac{n}{2}\right)^m e^{-2mt} (\|b_t\|_{2n}^2)^{n-m}(\|c_t\|_{2n}^2)^m. \]
Now, suppose $t\ge t_J(2,2n) = \frac{1}{2}\log n$.  Then $e^{-2mt} \le n^{-m}$.  Since $b,c\in\H(J,\left.\eps\right|_J)$ it follows from the inductive hypothesis that $\|b_t\|_{2n}\le\|b\|_2$ and $\|c_t\|_{2n}\le\|c\|_2$.  Thus,
\begin{eqnarray*}
\|a_t\|_{2n}^{2n} & \le & \|b\|_2^{2n} + \sum_{m=1}^n \binom{n}{m}\left(\frac{n}{2}\right)^m n^{-m} (\|b\|_2^2)^{n-m}(\|c\|_2^2)^m \\
                              & = & \left(\|b\|_2^2 + \frac{1}{2}\|c\|_2^2\right)^n,
\end{eqnarray*}
and from Equation \ref{a 2 norm}, this equals $\|a\|_2^{2n}$.  This proves the theorem. 
\end{proof}

\section{Speicher's Stochastic Interpolation}\label{Speicher}

In this final section, we consider creation and annihilation operators $\beta_j$, $\beta_j^\ast$ on $L^2(\CC,\te)$ which bear the same relation to the generators $x_j$ in $\CC$ as the creation and annihilation operators $c_q$, $c_q^\ast$ bear to the $q$-Gaussian variables $X_q\in\Gamm_q$. We use these operators, together with a non-commutative central limit theorem of Speicher, to approximate the $L^p(\H_q,\tau_q)$-norm by the norm on $L^p(\H,\te)$, and thus transfer Theorem \ref{eps strong hyp} from the context of the mixed spin holomorphic algebras to the arena of the $q$-holomorphic algebras, proving Theorem \ref{Main Theorem}.  All of the techniques in this section are analogs of Biane's ideas in \cite{Biane 1}.

\subsection{Creation and Annihilation operators on $L^2(\CC,\te)$}
Define operators $\beta_j$ on $L^2(\CC(I,\eps),\te)$ by
\[
\beta_j(x_A) = \begin{cases} x_jx_A,&\text{ if }j\notin A \\
                                        0,&\text{ if }j\in A
                           \end{cases}.
\]
One may readily verify that the adjoint of $\beta_j$ is given by
\[
\beta_j^\ast(x_A) = \begin{cases} 0,&\text{ if }j\notin A \\
                                                 x_jx_A,&\text{ if }j\in A \\
                                   \end{cases}.
\]
In the case $\eps(i,j) = 1$ for $i\ne j$, these are the \emph{B\'eb\'e Fock} operators on the toy Fock space of \cite{Meyer}.  In general, $\beta_j$ and $\beta_j^\ast$ mimic the creation and annihilation operators.  It is easy to see from dimension considerations that the $\ast$-algebra they generate is all of $\B(L^2(\CC(I,\eps),\te))$. We also have
\begin{equation}\label{beta + beta* = x}
\beta_j+\beta_j^\ast = x_j,
\end{equation}
as a left-multiplication operator on $\CC(I,\eps)$.  One can readily compute that these operators $\eps$-commute -- i.e.\ $\beta_i \beta_j = \sigma(i,j) \beta_j \beta_j$ if $i\ne j$.  They also satisfy the $\eps$-relations
\[ \beta_i^\ast \beta_j -\eps(i,j)\beta_j\beta_i^\ast = \delta_{ij}, \]
just like the operators $\z_j\in\CC(I\times\{0,1\},\eps)$.  In fact, the map $\z_j\mapsto\beta_j$ induces a $\ast$-isomorphism from $\CC(I\times\{0,1\},\eps)$ onto $\B(L^2(\CC(I,\eps)))$.  (In the case $\eps\equiv-1$, this reduces to the well known isomorphism from the complex Clifford algebra $\CC_{2n}$ onto the full matrix algebra $M_{2^n}(\C)$.)  Beware, however: this isomorphism does not send $\te$ to the normalized trace $tr$ on $\B(L^2(\CC(I,\eps)))$, as pointed out in Section \ref{holomorphic algebra}.

\medskip

The operators $\beta_j$, $\beta_j^\ast$ demonstrate concretely that the pure state $\beta\mapsto (\beta 1,1)_\eps$, the extension of $\te$ to $L^2(\CC(I,\eps),\te)$, is {\it not} tracial.  Indeed, it is easy to calculate that $(\beta_j\beta_j^\ast1,1)_\eps = 0$ while $(\beta_j^\ast\beta_j1,1)_\eps = 1$.  These are, however, the same covariance relations that the operators $c_q$ and $c_q^\ast$ satisfy with respect to the pure state $A\mapsto(A\vac,\vac)_q$ on $\B(\F_q)$.  It is additionally true that $(\beta_j1,1)_\eps = (\beta_j^\ast1,1)_\eps=0$, also in line with the operators $c_q$ and $c_q^\ast$.

\medskip

The following lemma shows that the state $(\cdot1,1)_\eps$ factors over naturally ordered products of the operators $\beta_j$ and $\beta_j^\ast$. It is proved in \cite{Biane 1}.
\begin{lemma}\label{factorize}
For each $j\in I$, let $\alpha_j$ be in the $\ast$-algebra generated by $\beta_j$.  Let $j_1,\ldots,j_s$ be $s$ distinct elements in $I$.  then
\[ (\alpha_{j_1}\cdots\alpha_{j_s}1,1)_\eps = (\alpha_{j_1}1,1)_\eps\cdots(\alpha_{j_s}1,1)_\eps. \]
\end{lemma}

\subsection{Speicher's central limit theorem}\label{clt}
Fix $q\in[-1,1]$.  We consider the family of random matrices ${\mathfrak S}_q$, consisting of all those infinite symmetric random matrices
$\eps\colon\N^\ast\times\N^\ast\to\{-1,1\}$ constantly $-1$ on the diagonal, for which $\{\eps(i,j)\,;\,i < j\}$ are i.i.d.\ with $\P(\eps = 1) = (1+q)/2$.
Note, then, $\P(\eps=-1) = (1-q)/2$, and so
\[ \E(\sigma(i,j)) = \frac{1+q}{2}\cdot1 + \frac{1-q}{2}\cdot-1 = q. \]
This family of random matrices features prominently in the main theorem of \cite{Speicher}, which we will use to prove Theorem \ref{Main Theorem}.

\medskip

Let $I_n$ denote the set $\{1,\ldots,n\}$, and let $\sigma\in{\mathfrak S}_q$.  For convenience, we denote the algebra $\CC(I_n\times\{0,1\},\left.\eps\right|_{I_n})$ as $\CC(n,\eps)$.  The creation operators  on $\CC(n,\eps)$ are labeled by pairs $(j,\zeta)$ where $j\in I_n$ and $\zeta\in\{0,1\}$;  to avoid confusion, we also index them as $\beta_{j,\zeta}^\eps$ to keep track of the dependence on $\eps$.  Let $d$ be a positive integer, and define new variables $\beta_{1,\zeta}^{\eps,n},\ldots,\beta_{d,\zeta}^{\eps,n}$, which act on $\CC(nd,\eps)$, by
\[ \beta^{\eps,n}_{k+1,\zeta} = \frac{1}{\sqrt{n}}\sum_{\ell = nk+1}^{n(k+1)} \beta_{\ell,\zeta}^\eps,\quad 0\le k \le d-1. \]
These  operators are constructed to approximate the operators $c_q$.  The intuition is: due to the expectation of the matrix $\eps\in{\mathfrak S}_q$, for large $n$ the $\beta^{\eps,n}_{j,\zeta}$ satisfy commutation relations close to the $q$-commutation relations of Equation \ref{q-commutation relations}.  Speicher's central limit theorem (Theorem 2 in \cite{Speicher}) makes this statement precise, but requires that the matrix of spins for the different variables have independent (upper triangular) entries.  In our case, since for each pair $i,j$ the entries $\eps((i,\zeta),(j,\zeta'))$ are the same for all choices of $\zeta,\zeta'\in\{0,1\}$, the matrix is only {\it block-independent} (with blocks of size $2\times 2$).  Nevertheless, as with the classical central limit theorem, a straighforward modification of Speicher's proof generalizes the theorem to this case.  We thus have the following theorem.
\begin{theorem}\label{stochastic}
Let $e_1,\ldots,e_d$ be an orthonormal basis for $\R^d$.  Among the operators $c_q(e_j,e_k)$ on $\F_q(\R^d\oplus\R^d)$, denote $c_q(e_j,0)$ as $c^q_{j,0}$, and denote $c_q(0,e_j)$ as $c^q_{j,1}$. Let $Q$ be a polynomial in $4d$ non-commuting variables.  For almost every $\eps\in{\mathfrak S}_q$,
\[\begin{aligned}
\lim_{n\to\infty} (&Q(\beta_{1,0}^{\eps,n},\ldots,\beta_{d,1}^{\eps,n},(\beta_{1,0}^{\eps,n})^\ast,\ldots,(\beta_{d,1}^{\eps,n})^\ast)1,1)_\eps \\
&= (Q(c^q_{1,0},\ldots,c^q_{d,1},(c^q_{1,0})^\ast,\ldots,(c^q_{d,1})^\ast)\vac,\vac)_q.
\end{aligned} \]
\end{theorem}
\begin{proof} This follows from Speicher's central limit theorem.  The required covariance conditions for the operators $\beta^{\eps}_{j,\zeta}$ were verified above, and the factorization of naturally-ordered products is the content of Lemma \ref{factorize}. 
\end{proof}
An immediate corollary is that the moments of elements in $\H_q(\C^d)$ can be approximated by the corresponding elements in $\CC(nd,\eps)$.  To be precise: let $x^{\eps,n}_j = \beta^{\eps,n}_{j,0} + (\beta^{\eps,n}_{j,0})^\ast$ and let $y^{\eps,n}_j = \beta^{\eps,n}_{j,1} + (\beta^{\eps,n}_{j,1})^\ast$.  By
Equation \ref{beta + beta* = x},
 \[ x^{\eps,n}_{j+1} = \frac{1}{\sqrt{n}}\sum_{\ell = nj+1}^{n(j+1)} x^\eps_\ell, \quad y^{\eps,n}_{j+1} = \frac{1}{\sqrt{n}}\sum_{\ell = nj+1}^{n(j+1)} y^\eps_\ell. \]
Let $z^{\eps,n}_j = 2^{-1/2}(x^{\eps,n}_j + iy^{\eps,n}_j)$, which is in $\H(I_{nd},\left.\eps\right|_{I_{nd}})$.
\begin{proposition}\label{stochastic 2}
Denote $Z_q(e_j)$ as $Z_j^q$. Let $r$ be an even integer, and let $P$ be a polynomial in $d$ non-commuting variables. For almost every
$\eps\in{\mathfrak S}_q$,
\[ \lim_{n\to\infty}\|P(z^{\eps,n}_1,\ldots,z^{\eps,n}_d)\|_{L^r(\H,\te)} = \|P(Z^q_1,\ldots,Z^q_d)\|_{L^r(\H_q,\tau_q)}. \]
\end{proposition}
\begin{proof}  Let $Q$ be the polynomial in $4d$ non-commuting variables defined by
\[  Q(\beta_{1,0}^{\eps,n},\ldots,\beta_{d,1}^{\eps,n},(\beta_{1,0}^{\eps,n})^\ast,\ldots,(\beta_{d,1}^{\eps,n})^\ast)
  = P(z^{\eps,n}_1,\ldots,z^{\eps,n}_d)^\ast P(z^{\eps,n}_1,\ldots,z^{\eps,n}_d).
\]
Such a polynomial exists because the variable $z^{\eps,n}_j$ is a (linear) polynomial in $\beta_{j,0}^{\eps,n}$, $\beta_{j,1}^{\eps,n}$, and their adjoints.  By definition, the same polynomial yields
\[ Q(c^q_{1,0},\ldots,c^q_{d,1},(c^q_{1,0})^\ast,\ldots,(c^q_{d,1})^\ast) = P(Z^q_1,\ldots,Z^q_d)^\ast P(Z^q_1,\ldots,Z^q_d). \]
Applying Theorem \ref{stochastic} to the polynomial $Q^m$, we have
\[ \lim_{n\to\infty} \te|P(z^{\eps,n}_1,\ldots,z^{\eps,n}_d)|^{2m}
= \tau_q|P(Z^q_1,\ldots,Z^q_d)|^{2m} \quad a.s.[\eps] \]
where we have used the fact that $(\cdot 1,1)_\eps$ reduces to $\te$ when applied to elements of $\CC(nd,\eps)$. 
\end{proof}

\medskip

We will also need to know that the semigroup $e^{-tN_\eps}$ approximates $e^{-tN_q}$.

\begin{proposition}\label{approx 2}
Let $r$ be an even integer, and let $P$ be a polynomial in $d$ non-commuting variables.  For $t > 0$, and for almost every $\eps\in{\mathfrak S}_q$,
\[ \lim_{n\to\infty} \|e^{-tN_\eps}P(z^{n,\eps}_1,\ldots,z^{n,\eps}_d)\|_r = \|e^{-tN_q}P(Z^q_1,\ldots,Z^q_d)\|_r. \]
\end{proposition}
\begin{proof} We can expand $P(Z^q_1,\ldots,Z^q_d)$ as a linear combination of monomials $Z^q_{i_1}\cdots Z^q_{i_\ell}$.  Each such monomial is an eigenvector of $e^{-tN_q}$ with eigenvalue $e^{-\ell t}$.  So it is easy to see that there is a unique polynomial
$P_t$ such that
\[ P_t(Z^q_1,\ldots,Z^q_d) = e^{-tN_q}P(Z^q_1,\ldots,Z^q_d). \]
Now, consider the polynomials $z^{\eps,n}_{i_1}\cdots z^{\eps,n}_{i_\ell}$.  Since $z^{\eps,n}_i$ is a linear combination of $z^\eps_1,\ldots,z^\eps_{nd}$, this polynomial may be expanded as a linear combination of monomials $z^\eps_{j_1}\cdots z^\eps_{j_\ell}$ with $1\le j_1,\ldots,j_\ell \le nd$.  From Equation
\ref{z commute}, if any two indices are equal, then $z^\eps_{j_1}\cdots z^\eps_{j_\ell}=0$; otherwise, it is of degree $\ell$.  Hence
\[ e^{-tN_\eps}(z^\eps_{j_1}\cdots z^\eps_{j_\ell}) = e^{-\ell t}(z^\eps_{j_1}\cdots z^\eps_{j_\ell}). \]
It follows that $e^{-tN_\eps}(z^{\eps,n}_{i_1}\cdots z^{\eps,n}_{i_\ell}) = e^{-\ell t}(z^\eps_{j_1}\cdots z^\eps_{j_\ell})$.  Thus, we see that
\[ P_t(z^{\eps,n}_1,\ldots,z^{\eps,n}_d) = e^{-tN_\eps}P(z^{\eps,n}_1,\ldots,z^{\eps,n}_d). \]
The theorem now follows by applying Proposition \ref{stochastic 2} to the polynomial $P_t$.  
\end{proof}
It should be noted that this elementary argument {\it fails} in the full algebra $\CC(nd,\eps)$; for example, $(x^\eps_1)^2 = 1$ is of degree $0$,
while $X_q(e_1)^2$ is of degree $2$ if $q>-1$.  The relevant statement is still true in that case, but a much more delicate argument (which can
be found in \cite{Biane 1}) is necessary to prove it.

\medskip

We now conclude with the end of the proof of Theorem \ref{Main Theorem}.
\begin{proof}[Proof of Theorem \ref{Main Theorem}]  First note that the sharpness of the Janson time $t_J(p,r)$ for any $p,r>0$ can be confirmed
by an argument identical to the one in the proof of Theorem \ref{eps strong hyp}.  For sufficiency, by standard arguments it is enough to prove the theorem for the finite dimensional Hilbert space $\HH = \R^d$, and moreover it suffices to prove it for elements $f\in L^2(\H_q,\tau_q)$ that are polynomials $f = P(Z^q_1,\ldots,Z^q_d)$ of the generators.  Let $r$ be an even integer, and let $t\ge t_J(2,r)$.  By Proposition \ref{approx 2},
\[ \|e^{-tN_q}f\|_r = \lim_{n\to\infty} \|e^{-tN_\eps}P(z^{\eps,n}_1,\ldots,z^{\eps,n}_d)\|_r \quad a.s.[\eps]. \]
By Theorem \ref{eps strong hyp} applied to the algebra $\H(I_{nd},\left.\eps\right|_{I_{nd}})$,
\[ \|e^{-tN_\eps}P(z^{\eps,n}_1,\ldots,z^{\eps,n}_d)\|_r \le \|P(z^{\eps,n}_1,\ldots,z^{\eps,n}_d)\|_2. \]
Finally, applying Proposition \ref{stochastic 2}, we have
\[ \lim_{n\to\infty}  \|P(z^{\eps,n}_1,\ldots,z^{\eps,n}_d)\|_2 = \| f \|_2 \quad a.s.[\eps]. \]
This completes the proof.

\end{proof}

\bigskip
\bigskip

\noindent {\bf Acknowledgement.} I would like to thank Len Gross for much insight, and for suggesting a thesis problem which led to this work.  I would also like to thank Philippe Biane, Claus Koestler, and Roland Speicher for useful conversations.


\bigskip

\begin{thebibliography}{99}

\bibitem[B1]{Biane 1} Biane, P.: \emph{Free hypercontractivity.} Commun. Math. Phys. {\bf 184}, 457-474 (1997)

\bibitem[B2]{Biane 2} Biane, P.: \emph{Segal-Bargmann transform, functional calculus on matrix spaces and the theory of semi-circular and circular systems.} J. Funct. Anal. {\bf 144}, 232-286 (1997)

\bibitem[BSZ]{BSZ} Baez, Segal, Zhou.: \emph{Introduction to algebraic and constructive quantum field theory.} Princeton Series in Physics.  Princeton University Press, Princeton, 1992.

\bibitem[BSp]{BSp}  Bozejko, M., Speicher, R.: \emph{An example of a generalized Brownian motion.}  Comm. Math. Phys. {\bf 137}, 519--531 (1991)

\bibitem[BKS]{BKS} Bozejko, M., K\"ummerer, B., Speicher, R.: \emph{q-Gaussian processes: non-commutative and
classical aspects.} Commun. Math. Phys. {\bf 185}, 129-154 (1997)

\bibitem[CL]{Carlen Lieb} Carlen, E., Lieb, E.: \emph{Optimal hypercontractivity for Fermi fields and related non-commutative
integration inequalities.} Commun. Math. Phys. {\bf 155}, 26-46 (1993)

\bibitem[DN]{DN} Dykema, K.; Nica, A.: \emph{On the Fock representation of the $q$-commutation relations.}  J. Reine Angew. Math. {\bf 440}, 201-212 (1993)

\bibitem[G1]{Gross 1} Gross, L.: \emph{Existence and uniqueness of physical ground states.} J. Funct. Anal. {\bf 10}, 52-109 (1972)

\bibitem[G2]{Gross 2} Gross, L.: \emph{Logarithmic Sobolev inequalities.} Amer. J. Math. {\bf 97}, 1061-1083 (1975)

\bibitem[G3]{Gross 3} Gross, L.: \emph{Hypercontractivity and logarithmic Sobolev inequalities for the Clifford-Dirichlet form.}
Duke Math. J. {\bf 42}, 383-396 (1975)

\bibitem[G4]{Gross 4} Gross, L.: \emph{Hypercontractivity over complex manifolds.} Acta. Math. {\bf 182}, 159-206 (1999)

\bibitem[G5]{Gross 5} Gross, L.: \emph{Hypercontractivity, logarithmic Sobolev inequalities and applications: a survey of surveys.} To appear.

\bibitem[Gli]{Glimm} Glimm, J.: \emph{Boson fields with nonlinear self-interaction in two dimensions.} Commun. Math. Phys. {\bf 8}, 12-25 (1968)

\bibitem[J]{Janson} Janson, S.: \emph{On hypercontractivity for multipliers on orthogonal polynomials.} Ark. Math. {\bf 21}, 97-110 (1983)

\bibitem[JPR]{JPR} Janson, S.; Peetre, J.; Rochberg, R.: \emph{Hankel forms and the Fock space.}  Rev. Mat. Iberoamericana {\bf 3} no. 1, 61--138 (1987)

\bibitem[K]{Krolak} Krolak, I.: \emph{Contractivity properties of the Ornsetein-Uhlenbeck semigroup for
general commutation realtions.}  To appear in Mathematische Zeitschrift.

\bibitem[L]{Lindsay} Lindsay, J.M.: \emph{Gaussian hypercontractivity revisited.}  J. Funct. Anal. {\bf 92}, 313-324 (1990)

\bibitem[LM]{Lindsay Meyer} Lindsay, J.M., Meyer, P.A.: \emph{Fermionic hypercontractivity.} In: \emph{Quantum probability VII,} Singapore: Wold Scientific, 1992, pp. 211-220

\bibitem[M]{Meyer} Meyer, P.A.: \emph{Quantum Probability for Probabilists.} Second edition, \emph{Lecture Notes in Mathematics} Vol. {\bf 1538},
Berlin--Heidelberg--New York: Springer, 1995

\bibitem[N1]{Nelson 1} Nelson, E.: \emph{A quartic interaction in two dimensions.}  In:  \emph{Mathematical Theory of Elementary
Particles,} M.I.T. Press, 1965, pp. 69-73

\bibitem[N2]{Nelson 2} Nelson, E.: \emph{The free Markov field.} J. Funct. Anal. {\bf 12}, 211-227 (1973)

\bibitem[PX]{PX} Pisier, G., Xu, Q.: \emph{Non-commutative $L^p$-spaces.} In: \emph{Handbook of the geometry of Banach spaces, Vol. 2,} North-Holland, Amsterdam, 2003, pp. 1459-1517

\bibitem[R]{Ricard} Ricard, E.: \emph{Factoriality of $q$-Gaussian von Neumann algebras.} To appear in Commun. Math. Phys.

\bibitem[S]{Speicher} Speicher, R.: \emph{A non-commutative central limit theorem.} Math Zeit. {\bf 209}, 55-66 (1992)

\bibitem[Se1]{Segal 1} Segal, I.: \emph{Tensor algebras over Hilbert spaces, II.} Ann. of Math. {\bf 63}, 160-175 (1956)

\bibitem[Se2]{Segal 2} Segal, I.: \emph{Construction of non-linear local quantum processes, I.} Ann. of Math. {\bf 92}, 462-481 (1970)

\bibitem[\'Sn]{Sniady} \'Sniady, P.: \emph{Gaussian random matrix models for $q$-deformed Gaussian variables.}  Commun. Math. Phys., {\bf 216}, 515-537 (2001)

\bibitem[V]{Voiculescu} Voiculescu, D.V.: \emph{Symmetries of some reduced free product $C^\ast$ algebras.} In: \emph{Operator Algebras and their Connection with Topology and Ergodic Theory, Lecture Notes in Mathematics,} Vol. {\bf 1132}, Berlin-Heidelberg-New York: Springer, 1985, pp. 566-588

\end{thebibliography}
\end{document}